\documentstyle[amssymb]{article}

\def\z{\zeta}
\def\Z{{\Bbb Z}}
\def\C{{\Bbb C}}
\def\Q{{\Bbb Q}}
\def\N{{\Bbb N}}

\newcommand\shf{
\setlength{\unitlength}{.4pt}
\begin{picture}(40,20)
\put(10,2){\line(1,0){20}}
\put(10,2){\line(0,1){10}}
\put(20,2){\line(0,1){10}}
\put(30,2){\line(0,1){10}}
\end{picture}
}

\date{}

\def\A{{\cal A}}

\def\proof{ \medskip
\noindent
{\em Proof.\ \ }}

\def\eop {\hfill $\Box$ \par\medskip}

\newtheorem{Def}{Definition}
\newtheorem{Lem}{Lemma}
\newtheorem{Thm}{Theorem}
\newtheorem{Cor}{Corollary}
\newtheorem{Conj}{Conjecture}
\newtheorem{Prop}{Proposition}

\begin{document}

\title{\bf Combinatorial Aspects\\ 
of Multiple Zeta Values}
\maketitle
\begin{center}

\vspace{2\baselineskip}

{\sc Jonathan M.~Borwein\footnote{Research supported by NSERC and the Shrum
Endowment of Simon Fraser University.}}
\par
{\em CECM, Department of Mathematics and Statistics, Simon Fraser University,
Burnaby, B.C., V5A$\;$1S6, Canada (e-mail: jborwein@cecm.sfu.ca)} 
\par \vspace{\baselineskip}
{\sc David M.~Bradley\footnote{Work done while the author was recipient of
the 
NSERC Postdoctoral Fellowship.}}  \par
{\em Department of Mathematics and Statistics, Dalhousie University,
Halifax, N.S., B3H$\;$3J5, Canada (e-mail: bradley@mscs.dal.ca)} 
\par \vspace{\baselineskip}
{\sc David J.~Broadhurst} \par
{\em Physics Department, Open University, Milton Keynes, MK7$\;$6AA, UK
     (e-mail: D.Broadhurst@open.ac.uk)}
\par \vspace{\baselineskip}
{\sc Petr Lison\v ek\footnote{Industrial Postdoctoral Fellow of PIms 
(The Pacific Institute for the Mathematical Sciences).\\
\ \\
{\em AMS (1991) subject classification:} Primary 05A19, 11M99, 68R15,
Secondary 11Y99.}}\par
{\em CECM, Department of Mathematics and Statistics, Simon Fraser University,
Burnaby, B.C., V5A$\;$1S6, Canada (e-mail: lisonek@cecm.sfu.ca)}

\vspace{1.5cm}
Submitted: July 2, 1998; Accepted: August 1, 1998.

\end{center}

\newpage

\begin{abstract}
Multiple zeta values (MZVs,
also called Euler sums or multiple harmonic series) 
are nested generalizations of the classical Riemann zeta function
evaluated at integer values.
The fact that an integral representation of MZVs
obeys a shuffle product rule allows the possibility of a~combinatorial
approach to them.  Using this approach we prove a~longstanding conjecture 
of Don Zagier about MZVs with certain repeated arguments.  
We also prove a~similar cyclic sum identity.  Finally,
we present extensive computational evidence supporting an infinite
family of conjectured 
MZV identities that simultaneously generalize the Zagier identity.
\end{abstract}

\section{Introduction}
\label{sec-intro}

In this paper,
we continue our study of {\em multiple zeta values (MZVs),} sometimes
also called {\em Euler sums} or {\em Zagier sums,} defined by
\[
\zeta(s_1,\ldots,s_k)
:=
\sum_{n_1>n_2>\ldots>n_k>0} \quad\prod_{j=1}^k n_j^{-s_j}
\]
with $s_j\in\Z^+$ and $s_1>1$ to ensure the convergence.
The integer $k$ is called the {\em depth} of the sum $\zeta(s_1,\ldots,s_k)$.

MZVs can be generalized in many ways. In particular,
they are instances of  {\em multidimensional polylogarithms}
\cite{BBBLa}.  Such sums have recently attracted much attention, in part since
there are many fascinating identities among them. 
The  applications of MZVs involve some unexpected fields,
such as high energy physics and knot theory---see \cite{BBBLa} for
a~list of references.

Hoffman in his study \cite{H} of the $\ast$-product of MZVs (which 
we call the ``stuffle'' product in \cite{BBBLa})
 distinguishes  between
``algebraic'' and ``non-al\-ge\-bra\-ic'' relations
among MZVs---the latter ones involve a~limiting process in 
some essential way.
In the same spirit  we note that
some non-trivial MZV identities are consequences
of discrete (combinatorial) relationships involving the shuffle product.
  Hints that this may be the case
include the occurrence of binomial coefficients
(e.g., in (\ref{eul-zs-zt})). 
In the present paper we follow the combinatorial approach
by exploring the combinatorial content of the shuffle product
rule (\ref{weight-length-shuffle}) 
for the integral representation \cite{BBBLa} of MZVs.

In Section~\ref{sec-fact-id} we list some factorial identities
on which we base our later results.  In Section~\ref{sec-shuf-alg} 
we introduce the shuffle algebra and in Section~\ref{sec-sh-id}
we prove some combinatorial identities holding in this algebra.  
The relevance of the shuffle algebra for studying MZVs 
originates in the iterated integral representation of MZVs
which we briefly recall in Section~\ref{sec-int-rep}.
In Section~\ref{sec-proofs} we use shuffle
identities to prove the longstanding conjecture 
of Don Zagier \cite{Zag,BBB,BBBLa}:
\[   
\z(\{3,1\}^n) = \frac{2\pi^{4n}}{(4n+2)!}
\]
(where the notation $\{ X\}^n$ indicates $n$ successive instances
of the integer sequence $X$), as well as the similar ``dressed with~2'' 
identity:
\[
\sum_{\vec s} \zeta(\vec s) = {{\pi^{4n+2}}\over{(4n+3)!}},
\]
where $\vec s$ runs over all $2n+1$ possible insertions of the number 2
in the string $\{3,1\}^n$.  Finally, in Section~\ref{sec-conj}
we present extensive numerical
evidence for our new conjecture,
which in a~rotationally symmetric way generalizes (by insertions
of groups of 2's) the Zagier identity. 
For an illustration, one very simple instance of our conjecture reads
\[
\z(3,2,2,1,2)+\z(2,2,3,2,1)+\z(2,3,1,2,2)={{\pi^{10}}\over{11!}}.
\]

\section{Factorial Identities}
\label{sec-fact-id}

In the main part of the paper we will require the following
identities.  The proofs
are easy by any of several methods 
(generating functions, WZ theory, etc.); therefore we skip them.

\begin{Lem}
\label{L1}
For any non-negative integer $n$ we have
\[
\sum_{r=-n}^n\frac{(-1)^r}{(2n+2r+1)!(2n-2r+1)!} 
=
\frac{2^{2n+1}}{(4n+2)!}.
\]
\end{Lem}

\begin{Lem}
\label{L3}
For any non-negative integer $n$ we have
\begin{equation}
\label{sum-L3}
\sum_{r=0}^n \frac{(-1)^r (2r+1)}{(2n+1-2r)!(2n+3+2r)!} = \frac{4^n}{(4n+3)!}.
\end{equation}
\end{Lem}

\begin{Lem}
\label{L2}
For any non-negative integer $n$ we have
\begin{eqnarray}
\label{sum-L2}
\sum_{r=0}^{n} (-1)^r(2r+1){{2n+1}\choose{n-r}}
=
\left\{
\begin{array}{ll}
1 & \mbox{if $n=0$}\\
0 & \mbox{if $n>0$.}
\end{array}
\right.
\end{eqnarray}
\end{Lem}

\section{The Shuffle Algebra}
\label{sec-shuf-alg}

Let $\A$ denote a~finite {\em alphabet} (set of letters).
By a~{\em word} on the alphabet $\A$ we mean a~(possibly empty)
sequence of letters from~$\A$.
By $\A^*$ we denote the set of all words on the alphabet $\A$.
For $w\in\A^*$, let $w^k$ denote the sequence of $k$ 
consecutive occurrences of~$w$.
A~{\em polynomial} on $\A$ over $\Q$ is a~rational linear combination 
of words on~$\A$. The set of all such polynomials is denoted by
$\Q\langle \A\rangle$.

On $\Q\langle \A\rangle$ we introduce the binary operation $\shf$
{\em (``shuffle product''),} which is defined, for any $u,v\in \A^*$ 
($u=x_1\ldots x_n$ and $v=x_{n+1} \ldots x_{n+m}$,  
$x_k\in \A$ for $1\le k\le n+m$) by
\begin{equation}
\label{shf-rule}
u\shf v  := \sum  x_{\sigma(1)}x_{\sigma(2)}\ldots x_{\sigma(n+m)},
   \label{shuff}
\end{equation}
where the sum is over all ${n+m\choose n}$ permutations $\sigma$ of the set
$\{1,2,\ldots,n+m\}$ which satisfy $\sigma^{-1}(j)<\sigma^{-1}(k)$ for
all $1\le j<k\le n$ and $n+1\le j<k\le n+m$.  In other words, the sum
is over all words (counting multiplicity) of length $n+m$ 
in which the relative
orders of the letters  $x_1,\dots,x_n$ and 
$x_{n+1},\dots,x_{n+m}$ are preserved.  The definition (\ref{shf-rule})
extends linearly on the entire domain 
$\Q\langle \A\rangle\times \Q\langle \A\rangle$.

\medskip
{\em Example.}  Let $\A=\{A,B\}$.  In $\Q\langle \A\rangle$
we have
\[
2AB\shf (3BA-AB)= 12AB^2A+12BA^2B+2(AB)^2+6(BA)^2-8A^2B^2.
\]

\section{Identities Involving Shuffles $(AB)^p$ with $(AB)^q$}
\label{sec-sh-id}

Throughout the rest of the paper we assume that the alphabet $\A$
contains exactly two letters $A$ and $B$.

\begin{Def}
\label{def-spqi}
Let $p$, $q$ and $j$ be non-negative integers
subject to ${\rm min}(p,q)\ge j$.
Let $S_{p+q,j}$ denote the set 
of those words occurring in $(AB)^p\shf (AB)^q$
that contain the subword $A^2$ exactly $j$ times.
\end{Def}

Definition~\ref{def-spqi} is sound, since the
{\em set} $S_{p+q,j}$ is the same for any partition of the number
$p+q$ into two parts as long as both parts are greater than or equal to $j$.
This would of course not be true if we instead considered the full expansion
of $(AB)^p\shf (AB)^q$ (that is, counting the multiplicity of words):
see Proposition~\ref{ABp-ABq} in which we calculate these multiplicities
explicitly.

{\em Side remark.}
The set $S_{p+q,j}$ has cardinality ${{p+q}\choose{2j}}$.
Indeed, any word in $S_{p+q,j}$ can be considered to be partitioned
into $p+q$ consecutive blocks of length~2. 
Clearly, the locations of the subwords $A^2$ and $B^2$ are 
consistent with this partitioning.
Since there are $j$ blocks
containing $A^2$, they must be interlaced with another $j$~blocks
containing $B^2$, and the choice of the positions of
these $j+j=2j$ blocks together with the shuffle rule (\ref{shf-rule})
determines the rest of the word in question.  Therefore there are exactly
${{p+q}\choose{2j}}$ elements in $S_{p+q,j}$.

\begin{Def}
Let $p,q,j$ be as in Definition~\ref{def-spqi}.
By $T_{p+q,j}$ we will denote the sum of all words in $S_{p+q,j}$.
\end{Def}

\begin{Prop}
\label{ABp-ABq}
For any non-negative integers $p$ and $q$ we have 
\[
(AB)^p \shf (AB)^q
= \sum_{j=0}^{\min(p,q)} 4^j\cdot{{p+q-2j}\choose{p-j}}\cdot T_{p+q,j}.
\]
\end{Prop}

\proof Let $u$ be an arbitrary but fixed word from $S_{p+q,j}$. 
Let us see how many
times $u$ arises in $(AB)^p\shf (AB)^q$.
This is the same as counting in how many ways the letters of $u$
can be colored in two colors (blue letters coming from $(AB)^p$
and red letters coming from $(AB)^q$)
in a coloring that is consistent with the shuffle rule (\ref{shf-rule}).

There are $p+q$ $A$'s in $u$, of which $2j$ $A$'s are contained 
in factors $A^2$
and $p+q-2j$ $A$'s
are surrounded by $B$'s from both sides (or possibly from
one side if we are looking at the leading $A$).
Of the latter $p+q-2j$ ``single'' $A$'s, $p-j$ are colored blue.
Thus the coloring of the single $A$'s
contributes a factor of ${{p+q-2j}\choose{p-j}}$ to the multiplicity of $u$
in $(AB)^p\shf (AB)^q$.  There are exactly $j$ factors $A^2$
(and thus exactly $j$ factors $B^2$) in $u$, each of which can be colored
in two ways (blue-red or red-blue), thus contributing a factor of 
$2^j\cdot 2^j=4^j$ to the multiplicity of $u$
in $(AB)^p\shf (AB)^q$.  
What remains to do is to color the ``single''
$B$'s, whose coloring is now determined uniquely by the choices 
made so far, together with the shuffle rule (\ref{shf-rule}).  
\phantom{X} \eop

\begin{Cor}
\label{cor-zagier-shfl}
For any non-negative integer $n$ we have
\begin{equation}
\label{zagier-shfl}
\sum_{r=-n}^n (-1)^r \left[ (AB)^{n-r} \shf (AB)^{n+r}\right]   
=  4^n (A^2B^2)^n.
\end{equation}
\end{Cor}

\proof Using Proposition~\ref{ABp-ABq}, the left-hand side of 
(\ref{zagier-shfl})
is equal to
\[
\sum_{r=-n}^n (-1)^r 
\sum_{j=0}^{\min(n-r,n+r)} 4^j\cdot{{2n-2j}\choose{n-r-j}}\cdot T_{2n,j}
\]
which after reordering is 
\begin{equation}
\label{zagier-swap}
\sum_{j=0}^n 4^j\cdot T_{2n,j} 
\sum_{r=j-n}^{n-j} (-1)^r{{2n-2j}\choose{n-r-j}}.
\end{equation}
Putting $N:=n-j$ in the inner sum turns it into
\[
\sum_{r=-N}^N (-1)^r{{2N}\choose{N-r}}
\]
which is equal to 1 if $N=0$ (i.e.~$j=n$)
whereas for $N>0$ (i.e.~$j<n$) it is a disguise
of $(1-1)^{2N}\cdot (-1)^N$ which is 0.  
Thus, (\ref{zagier-swap}) is equal to $4^n\cdot T_{2n,n}$ 
which is indeed the right-hand side of (\ref{zagier-shfl}), and the
proof is finished.
\eop

\begin{Cor}
\label{M1-Shf}
For any non-negative integer $n$ we have
\begin{eqnarray}
\label{hyper-Z-special}
\sum_{r=0}^n (-1)^r (2r+1) 
\left[ (AB)^{n-r} \shf (AB)^{n+1+r}\right]  =  \nonumber \\
  4^n\cdot\left( \sum_{r=0}^n  (A^2B^2)^r AB (A^2B^2)^{n-r}
+  \sum_{r=1}^n  (A^2B^2)^{r-1} A^2BAB^2 (A^2B^2)^{n-r} \right).
\end{eqnarray}
\end{Cor}

\proof
First we show that, for any non-negative integer $n$, we have 
\begin{equation}
\label{eq-L3}
\sum_{r=0}^n (-1)^r (2r+1) \left[ (AB)^{n-r} \shf (AB)^{n+1+r} \right]
=
4^n  T_{2n+1,n}.
\end{equation}
As in the proof of Corollary~\ref{cor-zagier-shfl} we proceed in three steps:
(i) evaluating the shuffle products by Proposition~\ref{ABp-ABq},
(ii) swapping the sums, 
(iii) doing the inner sum.

Using Proposition~\ref{ABp-ABq}, the left-hand side of (\ref{eq-L3})
can be written as
\[
\sum_{r=0}^n (-1)^r (2r+1)
\sum_{j=0}^{n-r} 4^j\cdot{{2n+1-2j}\choose{n-r-j}}\cdot  T_{2n+1,j}
\]
which after reordering is equal to
\[
\sum_{j=0}^{n}4^j T_{2n+1,j}\cdot
\sum_{r=0}^{n-j}(-1)^r (2r+1){{2n+1-2j}\choose{n-r-j}}
\]
which by Lemma~\ref{L2} (with $n-j$ in the place of $n$) 
is equal to $4^n T_{2n+1,n}$.

Now $T_{2n+1,n}$ is the sum of words
arising in the shuffle
$(AB)^{n} \shf (AB)^{n+1}$
and containing $n$ factors $A^2$ and $n$ factors $B^2$.  Thus, there
is exactly one single $A$ and exactly one single $B$, which clearly
have to be adjacent, and thus forming a factor $AB$ or $BA$.
In the parentheses on the right-hand side of~(\ref{hyper-Z-special}),
the first summand accounts for 
those summands from $T_{2n+1,n}$ that contain $AB$,
while the second summand accounts for 
those summands from $T_{2n+1,n}$ that contain $BA$.
This completes the proof of~(\ref{hyper-Z-special}).
\phantom{X} \eop

\section{Integral Representation of MZVs}
\label{sec-int-rep}

Let us recall that we are working with the alphabet $\A=\{A,B\}$.
Throughout the rest of this paper we identify 
the letter $A$ with the differential form $dx/x$ 
and the letter $B$ with the differential form $dx/(1-x)$.

The MZV $\zeta(s_1,\ldots,s_k)$ admits 
the $(s_1+s_2+\cdots+s_k)$-dimensional
iterated integral representation
\begin{equation}
\label{int-rep}
\zeta(s_1,\ldots,s_k)
=\int_0^1 A^{s_1-1} B A^{s_2-1} B \cdots
A^{s_k-1}B,
\qquad s_1>1.\label{casir}
\end{equation}
The explicit observation that MZVs are values 
of iterated integrals is apparently
due to Maxim Kontsevich~\cite{Zag}.
Less formally, such representations go as far back as Euler.
The representation (\ref{int-rep})
 is a~very special instance of the iterated integral representation
of multidimensional polylogarithms \cite{BBBLa}---see there
for the exact definition of the 
iterated integral (\ref{int-rep}), which however 
is not critical for our purposes.

Indeed, the only property of iterated integrals that we use in this paper
is that their products obey
the ``shuffle rule,'' that is \cite{Ree,BBBLa}
\begin{equation}
\label{weight-length-shuffle}
\left(\int_0^1 U \right)\cdot \left(\int_0^1 V \right) 
= \int_0^1 \left( U\shf V \right)
\end{equation}
if we view the products of differential 1-forms in $U$ and $V$ as words
in the shuffle algebra (Section~\ref{sec-shuf-alg}).
Clearly, (\ref{weight-length-shuffle}) motivated
our interest in shuffle identities (Section~\ref{sec-sh-id}).

An intriguing aspect of (\ref{weight-length-shuffle}) 
is the bridge
between analytical (transcendental) and discrete  nature
of MZVs.  Although the present paper deals only with MZVs, the
ideas used here are applicable to more general nested sums 
(alternating sums, multidimensional polylogarithms \cite{BBBLa}) since, 
as already mentioned above, 
these sums admit integral representations which generalize (\ref{int-rep}).

{\em Example.}
We provide a~combinatorial derivation of Euler's decomposition formula 
($s,t\ge 2$)
\begin{eqnarray}
\label{eul-zs-zt}
\z(s)\z(t)&=&\sum_{j=1}^s {{s+t-j-1}\choose {s-j}}\z(s+t-j,j)\nonumber\\
          && +\sum_{j=1}^t {{s+t-j-1}\choose {t-j}}\z(s+t-j,j).
\end{eqnarray}
Let us consider the product $P:=A^{s-1}B\shf A^{t-1}B$.
Clearly, any term in $P$ must end with a $B$. The terms in $P$
in which the trailing $B$ comes from the $A^{s-1}B$ operand are accounted for
by
\begin{equation}
\label{terms-1}
\sum_{k=t}^{s+t-1} {{k-1}\choose{t-1}} A^{k-1}BA^{s+t-k-1}B,
\end{equation}
with the binomial coefficient counting the number of ways in which 
all $A$'s from the $A^{t-1}B$ operand can be inserted in the leading block
of $A$'s in the shuffled string.
Similarly, those terms in $P$
in which the trailing $B$ comes from the $A^{t-1}B$ operand are accounted for
by
\begin{equation}
\label{terms-2}
\sum_{k=s}^{s+t-1} {{k-1}\choose{s-1}} A^{k-1}BA^{s+t-k-1}B.
\end{equation}
Summing up (\ref{terms-1}) and (\ref{terms-2}), 
substituting $k:=s+t-j$ and using 
(\ref{weight-length-shuffle},\ref{int-rep}) 
gives (\ref{eul-zs-zt}).
\section{Proof of the Zagier Conjecture}
\label{sec-proofs}
{}From Section~\ref{sec-intro} we recall that, 
in the context of integer sequences,
we use the notation $\{ X\}^n$ to indicate $n\ge0$ successive instances
of the sequence $X$.

\begin{Thm}
{\bf (The Zagier Conjecture)}
For any positive integer $n$ we have
\begin{equation}
\label{ZagC}
   \z(\{3,1\}^n) = \frac{2\pi^{4n}}{(4n+2)!}.
\end{equation}
\end{Thm}

\proof 
Using (\ref{weight-length-shuffle},\ref{int-rep}), 
Corollary~\ref{cor-zagier-shfl} implies
\[
\sum_{r=-n}^n (-1)^r\z(\{2\}^{n-r})\z(\{2\}^{n+r}) = 4^n\z(\{3,1\}^n).
\]
Application of the evaluation
\begin{equation}
\label{eval-twor}
\z(\{2\}^r)=  {{\pi^{2r}}\over{(2r+1)!}},
\end{equation}
which was proven in \cite{H-Pac,BBB}, gives
\[
4^n\z(\{3,1\}^n)=\pi^{4n}\sum_{r=-n}^n {{(-1)^r}\over{(2n-2r+1)!(2n+2r+1)!}}
\]
which by Lemma~\ref{L1} is equivalent to
\[
4^n\z(\{3,1\}^n)=\pi^{4n}\frac{2^{2n+1}}{(4n+2)!}.
\]
After dividing the last equation by $4^n$ we get
(\ref{ZagC}). \phantom{X}
\eop
The first proof of (\ref{ZagC}) appears in \cite{BBBLa}.
It may be viewed as the first non-commutative extension
of Euler's evaluation of $\z(2n)$.

\begin{Thm}
Let $n$ be a positive integer,
and let $I$ denote the set of all $2n+1$ possible insertions of the number 2
in the string $\{3,1\}^n$.
 Then
\begin{equation}
\label{dressed-eq}
\sum_{\vec s\in I} \zeta(\vec s) = {{\pi^{4n+2}}\over{(4n+3)!}}.
\end{equation}
\end{Thm}

\proof 
Using (\ref{weight-length-shuffle},\ref{int-rep}), 
Corollary~\ref{M1-Shf} implies
\[
\sum_{r=0}^n(-1)^r(2r+1)\z(\{2\}^{n-r})\z(\{2\}^{n+1+r})
=4^n \sum_{\vec s\in I} \zeta(\vec s).
\]
Indeed, the first term in the parentheses on the right-hand side 
of~(\ref{hyper-Z-special}) translates to 
$\sum_{r=0}^n \z(\{3,1\}^r,2,\{3,1\}^{n-r})$
while the second term translates to\\
$\sum_{r=1}^n \z(\{3,1\}^{r-1},3,2,1,\{3,1\}^{n-r})$.
Application of (\ref{eval-twor}) gives
\[
4^n \sum_{\vec s\in I} \zeta(\vec s)
=\pi^{4n+2}\sum_{r=0}^n {{(-1)^r(2r+1)}\over{(2n-2r+1)!(2n+3+2r)!}}
\]
which by Lemma~\ref{L3} is equivalent to
\[
4^n \sum_{\vec s\in I} \zeta(\vec s)=\pi^{4n+2}{{4^n}\over{(4n+3)!}}.
\]
After dividing the last equation by $4^n$ we get
(\ref{dressed-eq}). \phantom{X}
\eop

\section{Conjectured Generalizations of the Zagier\\ 
Identity}
\label{sec-conj}

To notationally ease our generalization, we define
\begin{equation}
\label{Z-def}
Z(m_0,\ldots,m_{2n})
:=
\zeta(\{2\}^{m_0},3,\{2\}^{m_1},1,\{2\}^{m_2},
\ldots,3,\{2\}^{m_{2n-1}},1,\{2\}^{m_{2n}})\,,
\end{equation}
with $\{2\}^{m_j}$ inserted after the $j$-th element of the string
$\{3,1\}^n$.
For example, $Z(2,0,1)=\z(2,2,3,1,2)$.

\begin{Conj}
\label{c-c}
For any sequence $S=(m_0,\ldots,m_{2n})$
 of $2n+1$ non-negative integers,
we have
\begin{equation}
\label{cyclic-eq}
\sum_{j=0}^{2n}Z({\cal C}^j S)  =     % \stackrel{?}{=}
\frac{\pi^{4n+2M}}{(4n+2M+1)!}\,,
\end{equation}
where $M:=\sum_{i=0}^{2n} m_i$ 
and ${\cal C}$ is the cyclic permutation operator, that is,
\[
{\cal C}^j (m_0,\ldots,m_{2n}) := 
(m_{2n-j+1},\ldots,m_{2n},m_0,\ldots,m_{2n-j}).
\]
\end{Conj}

{\em Remark.} Taking into account (\ref{eval-twor}) 
we see that the right-hand side of~(\ref{cyclic-eq}) is equal
to $\z\left(\{2\}^{2n+M}\right)$.

In Section~\ref{sec-proofs} we proved 
(\ref{cyclic-eq})
for the cases $M=0$ and $M=1$.
For $n=0$, (\ref{cyclic-eq}) trivially reduces 
to the known evaluation (\ref{eval-twor}).
If all $m_i$'s are equal,  (\ref{cyclic-eq}) specializes to conjecture (18)
of \cite{BBB}.

Since MZV duality~\cite{Ka,Le} implies that 
\begin{equation}
\label{Z-dual}
Z(S)=Z(\widetilde{S}), 
\end{equation}
where
$\widetilde{S}:=(m_{2n},\ldots,m_0)$ 
is the reverse of $S$, the conjecture (\ref{cyclic-eq})
can be also reformulated as~a sum over all permutations in 
the dihedral group $D_{2n+1}$.  
In our formulation we sum over the cyclic group $C_{2n+1}$.

\subsection{Integer Relations}
\label{sec-ir}

An {\em integer relation}  \cite{BL} for a~vector of complex numbers 
$z\in\C^n$ 
is a~non-zero vector of integers~$a\in\Z^n$ such that 
\[
a_1z_1+\cdots+a_nz_n=0.
\]
Conjecture~\ref{c-c} was discovered numerically (via its special instances) 
using the PSLQ algorithm for discovering integer relations \cite{FBA}
and the fast method for numerical evaluation of MZVs using
the H\"older convolution \cite{BBBLa}.  
All cases of (\ref{cyclic-eq})
with depth $2n+M\le 13$ were checked numerically at the precision
of 2000 digits. 
This amounted to checking 747 such identities, even after excluding
the cases with $n=0$ or $M\le 1$, for which proofs have been known before
or are presented in this paper.

\subsection{Other Conjectured Identities}

For any two fixed integers $n,M\ge 0$, let us consider
the vector $V_{n,M}$
of values $Z(m_0,m_1,\ldots,m_{2n})$ defined by (\ref{Z-def}) 
and subject to: $m_i\in\Z$, $m_i\ge 0$ ($0\le i\le 2n$) and
$\sum_{i=0}^{2n} m_i=M$.
We assume that the entries of $V_{n,M}$
are listed in some arbitrary (but fixed) order, and that
of any two $Z$-terms related by the duality (\ref{Z-dual}),
exactly one is present in $V_{n,M}$, in order to exclude trivial
duplicates. Additionally, we append to $V_{n,M}$
the value $Z(2n+M):=\z\left(\{2\}^{2n+M}\right)$.

If we restrict our attention to the {\em putative} identities
of the form (\ref{cyclic-eq}), then 
the number of (linearly independent) relations of this type
can be computed via P\'olya Theory (see, e.g.,~\cite{Ke}) 
as the number of orbits in the action
of the dihedral group $D_{2n+1}$ on the set of functions
$f:\{0,1,\ldots,2n\}\to\N$ subject to $\sum_{i=0}^{2n} f(i)=M$.
(Let us recall from Section~\ref{sec-ir} that we have verified
(\ref{cyclic-eq}) numerically in the range $2n+M\le 13$.)

On the other hand, integer relations for $V_{n,M}$ can be discovered
empirically using integer relation algorithms, regardless of whether
their structure is compatible with (\ref{cyclic-eq}) or not.
In Figure~\ref{fig1} we list, for some modest values of $n$ and $M$,
in lightface the number of (\ref{cyclic-eq})-type putative
relations for $V_{n,M}$,
and in boldface the number of relations for $V_{n,M}$ detected empirically
using the PSLQ algorithm \cite{FBA} using the numerical precision
of 5000 decimal places.
(In both cases we count the number of linearly independent relations.)
These values (as well as some others, not included in Figure~\ref{fig1})
suggest that the scheme (\ref{cyclic-eq}) exhaustively describes
all integer relations for $V_{n,M}$ in the cases 
when $n\le 1$ or $M\le 2$,
while in the remaining cases, additional relations were detected.

\newpage
\begin{figure}[h]
\centerline{
\begin{tabular}{|rr|r|r|r|}
\hline
 & $M$ & 1 & 2 & 3  \\
$n$ & & & & \\
\hline
1& & 1,{\bf 1} & 2,{\bf 2} & \phantom{0}3,{\bf \phantom{0}3} \\
\hline
2& & 1,{\bf 1}& 3,{\bf 3} & \phantom{0}5,{\bf \phantom{0}6} \\
\hline
3& & 1,{\bf 1}& 4,{\bf 4} & \phantom{0}8,{\bf 10} \\
\hline
4& & 1,{\bf 1}& 5,{\bf 5} & 12,{\bf 15} \\
\hline
\end{tabular}
}
\caption{Number of linearly independent integer relations for $V_{n,M}$.}
\label{fig1}
\end{figure}

{\em Example.}
Here is a~family of identities among $Z$-values (but not of type
(\ref{cyclic-eq})), for which we have extensive numerical evidence:

\begin{Conj}
For any non-negative integers $a_1,a_2,a_3,b_1,b_2$, we have
\begin{eqnarray*}
&&Z(a_1,b_1,a_2,b_2,a_3) + Z(a_2,b_1,a_3,b_2,a_1)  + Z(a_3,b_1,a_1,b_2,a_2)\\
&=&
Z(a_1,b_2,a_2,b_1,a_3) + Z(a_2,b_2,a_3,b_1,a_1)  + Z(a_3,b_2,a_1,b_1,a_2).
\end{eqnarray*}
\end{Conj}

\section*{Acknowledgment}

The authors gratefully acknowledge the support of the 
Canadian High Performance Computing Network (HPCnet), now merged with C3.ca.
Under this grant we started the development of 
{\em EZ-Face} (an abbreviation for Euler Zetas interFace), an on-line
calculator for Euler sums (by which we mean alternating MZVs),
 which is available 
for public use via the World Wide Web at the URL
\vskip.1in
\centerline{{\tt http://www.cecm.sfu.ca/projects/EZFace/}}


\vskip.1in



\begin{thebibliography}{99}


\bibitem{BBB}
J.~M.~Borwein, D.~M.~Bradley and D.~J.~Broadhurst,
``Evaluations of $k$-fold Euler/Zagier sums: A~compendium of results
for arbitrary $k$,''
{\em Electron.\ J.~Combin.}
{\bf 4} (1997), No.~2, \#R5.


\bibitem{BBBLa}
J.~M.~Borwein, D.~M.~Bradley, D.~J.~Broadhurst and P.~Lison\v ek,
``Special values of multidimensional polylogarithms,''
submitted.

\bibitem{BL}
J.~M.~Borwein and P.~Lison\v ek,
``Applications of integer relation algorithms,''
submitted.

\bibitem{FBA}
H.~R.~P.~Ferguson, D.~H.~Bailey and S.~Arno,
``Analysis of PSLQ, an integer relation finding algorithm,'' 
{\em Math.~Comp.,} to appear.  

\bibitem{H-Pac}
M.~E.~Hoffman,
``Multiple harmonic series,''
{\em Pacific J.\ Math.} {\bf 152} (1992), 275--290.


\bibitem{H}
M.~E.~Hoffman, ``The algebra of multiple harmonic series,''
{\em J.~Algebra} {\bf 194} (1997), 477--495.



\bibitem{Ka}
C.~Kassel,
``Quantum Groups,'' Springer, New York, 1995.


\bibitem{Ke}
A.~Kerber, ``Algebraic Combinatorics via Finite Group Actions,''
Bibliographisches Institut, Mannheim, 1991. 

\bibitem{Le}
T.~Q.~T.~Le and J.\ Murakami, 
``Kontsevich's integral for the Homfly polynomial and relations
between values of multiple zeta functions,''
{\em Topology Appl.} {\bf 62} (1995), 193--206.


\bibitem{Ree}
R.~Ree,
``Lie elements and an algebra associated with shuffles,''
{\em  Annals of Math.}
{\bf 68} (1958), 210--220.


\bibitem{Zag}
D.~Zagier, ``Values of zeta functions and their applications,''
{\em First European Congress of Mathematics},
Vol.~II, Birkh\"auser, Boston, 1994, pp.~497--512.

\end{thebibliography}
\end{document}